\documentclass[11pt]{article}
\input epsf.tex
\usepackage{amssymb,latexsym,times}
\catcode`\@=11
\renewcommand\subsection{\@startsection{subsection}{2}{\z@}%
                                     {-3.25ex\@plus -1ex \@minus -.2ex}%
                                     {-0.01 mm}
                                     {\normalfont\large\bfseries}}

\catcode`\@=11
\renewcommand\subsubsection{\@startsection{subsubsection}{2}{\z@}%
                                     {-3.25ex\@plus -1ex \@minus -.2ex}%
                                     {-0.01 mm}
                                     {\normalfont\bfseries}}

\newtheorem{example}{Example}
\newtheorem{theorem}[example]{Theorem}

\def\resp{{\em resp.$\ $}}

\def\cqfd{\hfill $\Box$ \bigskip}
\def\adots{\mathinner{\mkern2mu\raise1pt\hbox{.}
\mkern3mu\raise4pt\hbox{.}\mkern1mu\raise7pt\hbox{.}}}
\def\<{\langle\,}
\def\>{\,\rangle}

\def\ie{{\em i.e. }}
\def\eg{{\em e.g. }}

\def\SG{\mathfrak S}

\def\l{\lambda}
\def\a{\alpha}

\def\b{\beta}
\def\ga{\gamma}
\def\i{\iota}
\def\N{{\mathbb N}}
\def\Z{{\mathbb Z}}
\def\C{{\mathbb C}}

\def\Q{{\mathbb Q}}

\def\F{{\cal F}}

\def\B{{\mathbf B}}

\def\g{\mathfrak g}

\def\Sl{\mathfrak{sl}}

\def\n{\mathfrak n}

\def\slchap{\widehat{\mathfrak{sl}}}

\def\H{\widehat H}

\def\L{\Lambda}

\def\<{\langle}
\def\>{\rangle}

\def\CC{{\cal C}}
\def\m{\mu}
\def\mm{{\mathbf m}}
\def\nn{{\bf n}}

\def\le{\leqslant}
\def\ge{\geqslant}

\def\Si{\Sigma}
\def\si{\sigma}

\def\ind{{\rm Ind}}

\def\ra{\rightarrow}

\newdimen\Squaresize \Squaresize=14pt
\newdimen\Thickness \Thickness=0.5pt
 
\def\Square#1{\hbox{\vrule width \Thickness
   \vbox to \Squaresize{\hrule height \Thickness\vss
      \hbox to \Squaresize{\hss#1\hss}
   \vss\hrule height\Thickness}
\unskip\vrule width \Thickness}
\kern-\Thickness}
 
\def\Vsquare#1{\vbox{\Square{$#1$}}\kern-\Thickness}

\def\young#1{
\vbox{\smallskip\offinterlineskip
\halign{&\Vsquare{##}\cr #1}}}
 
\title{\bf A Littlewood-Richardson rule for evaluation representations
of $U_q(\slchap_n)$}

\author{Bernard {\sc Leclerc}}

\date{}

\begin{document}
\maketitle

\vskip 1cm

\begin{abstract}\noindent
We give a combinatorial description of the composition
factors of the induction product
of two evaluation modules of the affine Iwahori-Hecke algebra
of type $GL_m$.
Using quantum affine Schur-Weyl duality, this yields a combinatorial description 
of the composition factors of the tensor product of two
evaluation modules of the quantum affine algebra $U_q(\slchap_n)$. 
\end{abstract}

\vskip 0.6cm

\section{Introduction} \label{SECT1}

\subsection{}\label{1.1}
Let $H_m$ denote the Iwahori-Hecke algebra of type $A_{m-1}$
over $\C(t)$.
This is a semisimple associative algebra isomorphic
to the group algebra $\C(t)[\SG_m]$ of the symmetric group.
Hence its simple modules $S(\l)$ are parametrized by the 
partitions $\l$ of~$m$.
Consider a decomposition $m=m_1+m_2$, and two partitions
$\l^{(1)}$ and $\l^{(2)}$ of $m_1$ and $m_2$, respectively.
Then we have a $H_{m_1}$-module $S(\l^{(1)})$
and a $H_{m_2}$-module $S(\l^{(2)})$, and we can form
the induced module
\[
S(\l^{(1)})\odot S(\l^{(2)}) :=
\ind_{H_{m_1}\otimes H_{m_2}}^{H_m}\left( S(\l^{(1)})\otimes S(\l^{(2)})\right).
\]
Here, $H_{m_1}\otimes H_{m_2}$ is identified to a subalgebra
of $H_m$ in the standard way.
Using again the isomorphism $H_m\cong\C(t)[\SG_m]$, we see
that the multiplicity of a simple $H_m$-module $S(\m)$
in $S(\l^{(1)})\odot S(\l^{(2)})$ is equal 
to the classical Littlewood-Richardson coefficient 
$c_{\l^{(1)}\l^{(2)}}^\mu$ (see \eg \cite{Mcd}).

\subsection{}
Let now $\H_m$ be the affine Iwahori-Hecke algebra over $\C(t)$
(see \ref{evaluation} below).
For each invertible $z\in\C(t)$ we have a surjective 
{\em evaluation homomorphism}
$\tau_z : \H_m \ra H_m$.
Pulling back the simple $H_m$-module $S(\l)$ via $\tau_z$ we obtain
a simple $\H_m$-module $S(\l;z)$ called an {\em evaluation module}.
In analogy with \ref{1.1}, given two invertible elements
$z_1$ and $z_2$ of $\C(t)$, we can then form the induced $\H_m$-module
\[
S(\l^{(1)};z_1)\odot S(\l^{(2)};z_2) :=
\ind_{\H_{m_1}\otimes \H_{m_2}}^{\H_m}
\left( S(\l^{(1)};z_1)\otimes S(\l^{(2)};z_2)\right).
\]
It turns out that if we fix $\l^{(1)},\l^{(2)}$ and vary
the spectral parameters $z_1, z_2$, this module is generically
irreducible, that is, it is simple except for a 
finite number of values of the ratio $z_1/z_2$.
In \cite[Theorem 36]{LNT} a combinatorial description of
these special values was given.

In this note we shall make this result more precise by
describing all the composition factors of 
$S(\l^{(1)};z_1)\odot S(\l^{(2)};z_2)$ 
at these critical values $z_1/z_2$.
We shall also prove that, in contrast with the classical
Littlewood-Richardson rule, all the composition factors
appear with multiplicity one.
The composition factors occuring in a product will be
described using the combinatorics of Lusztig's {\em symbols}, that is,
of certain two-row arrays introduced by Lusztig for parametrizing
the irreducible complex representations of the classical
reductive groups over finite fields \cite{Lu79,Lu84}. 

\subsection{}
We will derive our combinatorial formula from some explicit
calculations of canonical bases in level 2 
representations of the quantum algebra $U_v(\Sl_{n+1})$
performed in \cite{LM}. 
More precisely, by dualizing \cite[Theorem 3]{LM}, we get
a formula for the expansion of the product of two quantum
flag minors on the dual canonical basis of $U_v(\Sl_{n+1})$
(Theorem~\ref{th2}).
Using then Ariki's theorem as in \cite{LNT}, we obtain immediately
the above-mentioned Littlewood-Richardson rule for induction products
of two evaluation modules over affine Hecke algebras (Theorem~\ref{th1}).

Finally, by means of the quantum affine analogue of the Schur-Weyl
duality developed by Cherednik, Chari-Pressley and 
Ginzburg-Reshetikhin-Vasserot, we can deduce from this
rule a similar one for the tensor product of two evaluation
modules over the quantum affine algebra $U_q(\slchap_N)$. 

\section{Composition factors of induced $\H_m$-modules} \label{SECT2}

\subsection{}\label{evaluation}
Let $\H_m$ be the affine Hecke algebra of type $GL_m$ over $\C(t)$.
It has invertible generators $T_1,\ldots ,T_{m-1},y_1,\ldots ,y_m$
subject to the relations
\begin{eqnarray*}
&&T_iT_{i+1}T_i=T_{i+1}T_iT_{i+1},\hskip 1.1cm (1\le i\le m-2),\label{EQ_T1}\\
&&T_iT_j=T_jT_i,\hskip 2.9cm (\vert i-j\vert>1),\label{EQ_T2}\\
&&(T_i-t)(T_i+1)=0, \hskip 1.6cm(1\le i\le m-1),\label{EQ_T3}\\
&&y_iy_j=y_jy_i, \hskip 3.1cm (1\le i,j\le m),\label{EQ_YY}\\
&&y_jT_i=T_iy_j, \hskip 3cm (j \not= i,i+1),\\
&&T_iy_iT_i=t\,y_{i+1}, \hskip 2.5cm (1\le i\le m-1).\label{EQ_Y}
\end{eqnarray*}
The subalgebra $H_m$ generated by the $T_i$'s is the Iwahori-Hecke
algebra of type $A_{m-1}$.

For any invertible $z\in\C(t)$ we have a unique algebra homomorphism
$\tau_z : \H_m \ra H_m$ such that
\[
\tau_z(T_i)=T_i,\quad \tau_z(y_1) = z,\qquad (i=1,\ldots,m-1).
\]
This is called the {\em evaluation at $z$}.

We also have an algebra automorphism 
$\si_z : \H_m \ra \H_m$ such that
\[
\si_z(T_i)=T_i,\quad \si_z(y_i) = zy_i,\qquad (i=1,\ldots,m-1).
\]
This is called the {\em shift by $z$}.

\subsection{}\label{reduction}
As mentioned in the introduction, given two partitions
$\l^{(1)}$ and $\l^{(2)}$, the structure of 
the induced $\H_m$-module
\[
S(\l^{(1)};z_1)\odot S(\l^{(2)};z_2)
\]
depends essentially on the ratio $z_1/z_2$.
Indeed, by twisting this module with the shift automorphism
$\si_z$ we obtain the induced module
\[
S(\l^{(1)};zz_1)\odot S(\l^{(2)};zz_2).
\]
For example, it is known that if $z_1/z_2\not\in t^\Z$ then 
$S(\l^{(1)};z_1)\odot S(\l^{(2)};z_2)$ is irreducible.
Therefore, we can assume without loss of generality that
\begin{equation}\label{assumption}
z_i=t^{a_i},\quad a_i\in\Z,\quad a_i\ge\ell(\l^{(i)}),\qquad
(i=1,2),
\end{equation}
where as usual $\ell(\l)$ denotes the length of the partition
$\l$.
Since $S(\l^{(1)};z_1)\odot S(\l^{(2)};z_2)$ and
$S(\l^{(2)};z_2)\odot S(\l^{(1)};z_1)$ have the same composition
factors with the same multiplicities, we can also assume that 
$a_1\le a_2$.

\subsection{}
It will be convenient to write partitions in weakly {\em increasing}
order.
Given a partition $\l$ and an integer $a\ge\ell(\l)$ we can make
$\l$ into a non-decreasing sequence $(\l_1,\ldots ,\l_a)$ 
of length $a$ by setting $\l_j=0$ for $j=1,\ldots ,a-\ell(\l)$.
We can then associate to $(\l,a)$ the increasing sequence
\begin{equation}\label{beta}
\b=(\b_1,\ldots ,\b_a), \quad \b_j=j+\l_j.
\end{equation}
In this way, given $(\l^{(i)},a_i)\ (i=1,2)$ as in
\ref{reduction}, we obtain a {\em symbol}
\begin{equation}\label{symbol}
S=\pmatrix{\b^{(2)} \cr \b^{(1)}} = 
\pmatrix{\b^{(2)}_1,\ldots ,\b^{(2)}_{a_2} \cr 
\b^{(1)}_1,\ldots ,\b^{(1)}_{a_1}}.
\end{equation}
For example, the symbol attached to the pairs $((1,1,2),3)$ and $((2,3),5)$
is 
\[
S=\pmatrix{1 &2 &3 &6 &8\cr 2& 3& 5}.
\]

Conversely, given a symbol $S$, \ie a two-row array as in Eq.~(\ref{symbol})
with 
\[
1\le \b^{(i)}_1<\cdots<\b^{(i)}_{a_i}\quad (i=1,2),
\]
there
is a unique pair $(\l^{(i)},a_i)\ (i=1,2)$ whose symbol is~$S$.

\subsection{}
The symbol $S$ of Eq.~(\ref{symbol}) is said to be {\em standard}
if $\b^{(2)}_k \le \b^{(1)}_k$ for $k\le a_1$.
In \cite[\S 2.5]{LM} we have defined the {\em pairs} of a standard
symbol $S$, and the set $\CC(S)$ of all symbols $\Si$ obtained
from $S$ by permuting some of its pairs.
As shown in \cite[Lemma 9]{LM}, these notions are equivalent to the
notion of admissible involution of Lusztig \cite{Lu02}.
 
For the convenience of the reader we shall recall these
definitions.
Let $S={\b\choose\ga}$ be a standard symbol.
We define an injection $\psi : \ga \longrightarrow \b$ such
that $\psi(j)\le j$ for all $j\in\ga$. 
To do so it is enough to describe the subsets
\[
\ga^l = \{j\in\ga \mid \psi(j)=j-l\},\qquad (0\le l \le n).
\]
We set $\ga^0 = \ga \cap \b$ and for $l\ge 1$ we put
\[
\ga^l = \{j\in\ga - (\ga^0 \cup \cdots \cup \ga^{l-1}) \mid 
j-l\in \b-\psi(\ga^0 \cup \cdots \cup \ga^{l-1})\}.
\]
Observe that the standardness of $S$ implies that $\psi$ is
well-defined.
\begin{example}
{\rm
Take 
\[
S=\pmatrix{1 & 3 & 5 & 8 & 9 \cr
           3 & 6 & 7 & 10}\,.
\]
Then 
\[
\ga^0 = \{3\},\ \ga^1 = \{6,10\},\ \ga^2=\cdots =\ga^5=\emptyset,
\ \ga^6 = \{7\}.
\]
Hence 
\[
\psi(3)=3,\ \psi(6)=5,\ \psi(7)=1,\ \psi(10)=9.
\] 
}
\end{example}

The pairs $(j,\psi(j))$ with $\psi(j)\not = j$ (that is, 
$j\not \in \b\cap\ga$) will be called the pairs of $S$.
Given a standard symbol $S$ with $p$ pairs, we denote by
$\CC(S)$ the set of all symbols obtained from $S$ by 
permuting some pairs in $S$ and reordering the rows. 
We consider $S$ itself as an element of $\CC(S)$, hence
$\CC(S)$ has cardinality $2^p$.

\subsection{}\label{young}
Given a partition $\l$ and an integer $a$ we call {\em Young diagram
of $(\l,a)$} the Young diagram of $\l$ in which each cell $(i,j)$ 
is filled with the integer $i-j+a$.
For instance, if $\l=(2,3)$ and $a=5$ then the Young diagram of
$(\l,a)$ is 
\[
\young{4 & 5\cr 5 & 6 & 7\cr}
\]

The rows of the Young diagram of $(\l,a)$ yield a {\em multisegment}
\[
\mm(\l,a):=\sum_{1\le k\le a} [k,k+\l_k-1].
\]
This is a formal sum (or multiset) of intervals in $\Z$,
in which we discard the empty intervals corresponding
to the $k$'s with $\l_k=0$.
Thus, continuing with the same example, we have
\[
\mm((2,3),5)=[4,5]+[5,7].
\]

Similarly, we attach to a pair $(\l^{(i)},a_i)\ (i=1,2)$
or to its symbol $S$ the multisegment
\[
\mm(S) = \mm(\l^{(1)},a_1) + \mm(\l^{(2)},a_2).
\]

\subsection{}
To each multisegment 
\[
\mm:=\sum_k [\a_k,\b_k]
\]
is attached an irreducible $\H_m$-module $L_\mm$, where 
$m=\sum_k (\b_k+1-\a_k)$ (see \eg \cite[\S 2.1]{LNT}).

\subsection{}
Let us assume that the pair $(\l^{(i)},a_i)\ (i=1,2)$
satisfies the conditions of \ref{reduction}.
Let $\Si$ denote the symbol attached to this pair. 
We can now state:

\begin{theorem}\label{th1}
The composition factors of 
$S(\l^{(1)};t^{a_1})\odot S(\l^{(2)};t^{a_2})$
are the modules $L_{\mm(S)}$ where $S$ runs through
the set of standard symbols such that $\Si\in\CC(S)$.
Each of them occurs with multiplicity one.
\end{theorem}
Theorem~\ref{th1} will be deduced from Theorem~\ref{th2}
below.

\begin{example}\label{exa1}
{\rm
Let $(\l^{(1)},a_1)=((1,4),2)$ and
$(\l^{(2)},a_2)=((1,2,3),4)$.
The corresponding symbol is
\[
\Si=\pmatrix{1 &3 &5 &7\cr 2& 6}.
\]
The standard symbols $S$ such that $\Si\in\CC(S)$ are
\[
\pmatrix{1 &2 &5 &6\cr 3& 7},\quad 
\pmatrix{1 &2 &5 &7\cr 3& 6},\quad
\pmatrix{1 &3 &5 &6\cr 2& 7},\quad 
\pmatrix{1 &3 &5 &7\cr 2& 6}. 
\]
It follows that the composition factors of 
$S((1,4);t^2)\odot S((1,2,3);t^4)$
are the $L_\mm$ where $\mm$ is one the following
multisegments:
\[
\nn_1=[1,2]+[2,6]+[3,4]+[4,5],\quad
\nn_2=[1,2]+[2,5]+[3,4]+[4,6],\quad
\]
\[
\nn_3=[1,1]+[2,2]+[2,6]+[3,4]+[4,5],\quad
\nn_4=[1,1]+[2,2]+[2,5]+[3,4]+[4,6].
\] 
}
\end{example} 
\section{Canonical bases} \label{SECT3}

\subsection{}
Fix $n \ge 2$ and let $\g=\Sl_{n+1}$. 
We consider the quantum enveloping algebra $U_v(\g)$ over
$\Q(v)$ with Chevalley generators 
$e_j, f_j, t_j \ (1\le j \le n)$.
The simple roots and the fundamental weights are denoted by 
$\a_k$ and $\L_k \ (1\le k \le n)$ respectively.
The irreducible representation of $U_v(\g)$ with highest
weight $\L$ is denoted by $V(\L)$.
We denote by $U_v(\n)$ the subalgebra of $U_v(\g)$ generated
by $e_j \ (1\le j \le n)$.

\subsection{}
Let $\B$ (\resp $\B^*$) denote the canonical basis (\resp the dual canonical
basis) of $U_v(\n)$ (\cite{Lu90}, \cite{BZ}; see also \cite[\S 3]{LNT}).
The elements of $\B$ and $\B^*$ are naturally labelled by
the multisegments $\mm$ supported on $[1,n]$. 
We shall denote them by $b_\mm$ and $b^*_\mm$ respectively.

The vectors $b^*_\mm$ for which $\mm$ is of the form
\[
\mm = \mm(\l,a)
\]
for some partition $\l$ and some integer $a$ are called 
{\em quantum flag minors}. 
Indeed, by \cite{BZ}, they can be expressed as quantum
minors of a triangular matrix whose entries are iterated
brackets of the $e_i$'s (see \cite[\S 5.2]{LNT}).

\subsection{}
Let $(\l^{(i)},a_i)\ (i=1,2)$ be as in \ref{reduction}.
We also assume that the multisegments 
\[
\mm_i = \mm(\l^{(i)},a_i)\quad (i=1,2)
\] 
are supported on $[1,n]$.
Let $\Si$ be the symbol attached to the pair $(\l^{(i)},a_i)\ (i=1,2)$.
For a standard symbol $S$ such that
$\Si\in\CC(S)$ we denote by $n(S,\Si)$ the number of pairs
of $S$ which are permuted to get $\Si$. 
Finally, we denote by $N_j(\l,a)$ the number of cells of 
the Young diagram of $(\l,a)$ containing the integer~$j$.

\begin{theorem}\label{th2}
We have
\[
b^*_{\mm_1}b^*_{\mm_2} =
v^{-N_{a_1}(\l^{(2)},\,a_2)}
\sum_S v^{n(S,\,\Si)}\,b^*_{\mm(S)}
\]
where the sum runs through all standard symbols $S$ such that
$\Si\in\CC(S)$.
\end{theorem}

\begin{example}\label{exa2}
{\rm
We take $(\l^{(1)},a_1)$ and $(\l^{(2)},a_2)$ as in Example~\ref{exa1}.
Hence
\[
\mm_1 = [1,1]+[2,5],\qquad \mm_2=[2,2]+[3,4]+[4,6].
\]
Then $N_{a_1}(\l^{(2)},\,a_2)=N_2((1,2,3),4)=1$, and we obtain,
using the notation of Example~\ref{exa1},
\[
b^*_{\mm_1}\, b^*_{\mm_2} =
v^{-1}
(v^2\,b^*_{\nn_1} + v\,b^*_{\nn_2} + v\,b^*_{\nn_3} + b^*_{\nn_4}).
\]
}
\end{example}

\subsection{} {\it Proof of Theorem~\ref{th2}.}
Following \cite[\S 7.2]{LNT}, we will replace
calculations of products of elements of $\B^*$ by calculations
of dual canonical bases of finite-dimensional representations
of $U_v(\g)$.

\subsubsection{} \label{part1}
Let $U_v(\n^-)$ denote the subalgebra of $U_v(\g)$ generated by
the $f_i$'s, and let $x \mapsto x^\sharp$ denote the algebra
isomorphism from $U_v(\n)$ to $U_v(\n^-)$ defined by 
$e_i^\sharp = f_i\ (i=1,\ldots, n)$.
Let $\L$ be a dominant integral weight and let $u_\L$
be a highest weight vector of the irreducible module $V(\L)$.
Then the map $\pi_\L : x\mapsto x^\sharp u_\L$ projects the canonical
basis $\B$ of $U_v(\n)$ to the union of the canonical basis $\B(\L)$
of $V(\L)$ with the set $\{0\}$.
The dual map $\pi_\L^*$ gives an embedding of the dual canonical
basis $\B^*(\L)$ of $V(\L)\simeq V(\L)^*$ into the dual canonical
basis $\B^*$ of $U_v(\n)\simeq U_v(\n)^*$.

\subsubsection{} \label{part2}
In particular the subset of $\B^*$ obtained by embedding the 
bases $\B^*(\L_a)\ (1\le a\le n)$ of the fundamental representations
is precisely the subset of quantum flag minors.
It is well known that $V(\L_a)$ is a minuscule representation
whose bases $\B^*(\L_a)$ and $\B(\L_a)$ coincide.
Moreover the elements of these bases are naturally labelled by
the pairs $(\l,a)$ whose Young diagram (as defined in \ref{young})
contains only cells numbered by integers between $1$ and $n$. 
Denoting them by $b^*_{(\l,a)}$ we have 
\[
\pi^*_{\L_a}(b^*_{(\l,a)}) = b^*_{\mm(\l,a)}
\]
Equivalently, we can also label the elements of $\B^*(\L_a)$
by one-row symbols $\b$ as in Eq.~(\ref{beta}) with $\b_i\le n+1$.

\subsubsection{} \label{part3}
Similarly, the basis $\B^*(\L_{a_1})\otimes \B^*(\L_{a_2})$ is naturally
labelled by the set of symbols $S$ as in Eq.~(\ref{symbol})
with $\b^{(i)}_{a_i}\le n+1\ (i=1,2)$.
Using the theory of crystal bases \cite{K1,K2} one can see
that the basis $\B^*(\L_{a_1}+\L_{a_2})$ has a natural labelling
by the subset of standard symbols \cite[\S 2.3]{LM}.
Moreover, denoting by $b^*_S$ the element of $\B^*(\L_{a_1}+\L_{a_2})$
labelled by the standard symbol $S$ we have, using also the
notation of \ref{young},
\[
\pi^*_{\L_{a_1}+\L_{a_2}}(b^*_S) = b^*_{\mm(S)}.
\]

\subsubsection{} \label{part4}
Let $\i : V(\L_{a_1}+\L_{a_2}) \ra V(\L_{a_1})\otimes V(\L_{a_2})$
be the $U_v(\g)$-module embedding which maps $u_{\L_{a_1}+\L_{a_2}}$
to $u_{\L_{a_1}}\otimes u_{\L_{a_2}}$, and let 
$\i^* : V(\L_{a_1})\otimes V(\L_{a_2}) \ra V(\L_{a_1}+\L_{a_2})$
be its dual. 
Let $b^*_i\in \B^*(\L_{a_i})\ (i=1,2)$
and denote by $b^*_{\mm_i} = \pi^*_{\L_{a_i}}(b^*_i)\ (i=1,2)$
the corresponding quantum flag minors.
It is shown in \cite[\S 7.2.7]{LNT} that the image of 
$b^*_1\otimes b^*_2$ under the composition of maps
$\pi^*_{\L_{a_1}+\L_{a_2}} \circ \i^*$ coincides up to a power of $v$
with the product $b^*_{\mm_1}b^*_{\mm_2}$.
Hence to calculate the $\B^*$-expansion of $b^*_{\mm_1}b^*_{\mm_2}$
it is enough to calculate the matrix of the map $\i^*$ with
respect to the bases $\B^*(\L_{a_1})\otimes\B^*(\L_{a_1})$
and $\B^*(\L_{a_1}+\L_{a_2})$.

\subsubsection{} \label{part5}
The matrix of $\i$ with respect to the bases
$\B(\L_{a_1}+\L_{a_2})$ and $\B(\L_{a_1})\otimes\B(\L_{a_1})$
was calculated in \cite[Theorem 3]{LM} in terms of Lusztig's symbols.
Transposing this matrix we obtain the desired matrix of $\i^*$.
Using \ref{part2} and \ref{part3}, we then get the formula
of Theorem~\ref{th2}.
\cqfd

\subsection{} {\it Proof of Theorem~\ref{th1}.}  
By \cite[\S 3.7]{LNT} the multisegments $\mm$
indexing the composition factors $L_\mm$ of 
$S(\l^{(1)};t^{a_1})\odot S(\l^{(2)};t^{a_2})$
are those occuring in the right-hand side of the formula
of Theorem~\ref{th2}.
Moreover the composition multiplicities are obtained by
specializing $v$ to $1$ in the coefficients of this formula.
Hence they are all equal to $1$.
\cqfd

\section{Tensor products of $U_q(\slchap_N)$-modules} \label{SECT4}

\subsection{} 
Let $U_q(\slchap_N)$ be the quantized affine algebra of
type $A_{N-1}^{(1)}$ with parameter $q$ a square root of~$t$
(see for example \cite{CP} for the defining relations of $U_q(\slchap_N)$).
The quantum affine Schur-Weyl duality between $\H_m$ and
$U_q(\slchap_N)$ \cite{CP, Ch, GRV} gives a functor $\F_{m,N}$ from
the category of finite-dimensional
$\H_m$-modules to the category of level $0$ finite-dimensional
representations of $U_q(\slchap_N)$.
If $N\ge m$, $\F_{m,N}$ maps the simple modules of $\H_m$
to simple modules of $U_q(\slchap_N)$. 
However, the image of a non-zero simple $\H_m$-module may be the
zero $U_q(\slchap_N)$-module.
More precisely, the simple $\H_m$-module $L_\mm$ is mapped
to a non-zero simple $U_q(\slchap_N)$-module if and only if all the segments 
occuring in $\mm$ have length $\le N-1$. 
In this case the Drinfeld polynomials of $\F_{m,N}(L_\mm)$
are easily calculated from $\mm$ (see \cite{CP}). 

The functor $\F_{m,N}$ transforms induction product into
tensor product, that is, 
for $M_1$ in $\CC_{m_1}$ and $M_2$ in $\CC_{m_2}$ one has
\[
\F_{m_1+m_2,N}(M_1\odot M_2) =
\F_{m_1,N}(M_1)\otimes\F_{m_2,N}(M_2)\,.
\]

\subsection{}
The image under $\F_{m,N}$ of an evaluation module for $\H_m$
is an evaluation module for $U_q(\slchap_N)$, and all evaluation
modules of $U_q(\slchap_N)$ can be obtained in this way,
by varying $m \in \N^*$.

\subsection{}
By application of the Schur functor $\F_{m,N}$ to Theorem~\ref{th1}
we thus obtain a combinatorial description of all composition
factors of the tensor product of two evaluation modules of
$U_q(\slchap_N)$.

\begin{example}
{\rm
We continue Example~\ref{exa1} and Example~\ref{exa2}.
The image of the $\H_5$-module $L_{\mm_1}$ under
$\F_{5,N}$ is the evaluation module $V(\mm_1)$ 
of $U_q(\slchap_N)$ with Drinfeld polynomials
\begin{eqnarray*}
P_1(u)& =& u - q^{-2},\\
P_2(u)& =& P_3(u)\ \  =\ \ 1, \\
P_4(u)& =& u - q^{-7},\\
P_k(u)& =& 1,\quad (5\le k\le N-1).
\end{eqnarray*}
This is a non-zero module if and only if $N\ge 5$.
Similarly, the image of the $\H_6$-module $L_{\mm_2}$ under
$\F_{6,N}$ is the evaluation module $V(\mm_2)$ 
of $U_q(\slchap_N)$ with Drinfeld polynomials
\begin{eqnarray*}
P_1(u) &=& u - q^{-4},\\
P_2(u) &=& u - q^{-7}, \\
P_3(u) &=& u - q^{-10},\\
P_k(u) &=& 1,\quad (4\le k\le N-1).
\end{eqnarray*}
This is a non-zero module if and only if $N\ge 4$.
The images of the $\H_{11}$-modules 
$L_{\mm_1}, L_{\mm_2}, L_{\mm_3}$, $L_{\mm_4}$
under $\F_{11,N}$
are the modules $V(\nn_1), V(\nn_2), V(\nn_3), V(\nn_4)$
with respective Drinfeld polynomials
\begin{eqnarray*}
\qquad\quad\
P_1(u) &=& 1,\\
P_2(u) &=& (u-q^{-3})(u-q^{-7})(u-q^{-9}),\\
P_3(u) &=& P_4(u)\ \  =\ \ 1, \\
P_5(u) &=& u - q^{-8},\\
P_k(u) &=& 1,\quad (6\le k\le N-1);
\end{eqnarray*} 
\begin{eqnarray*}
P_1(u) &=& 1,\\
P_2(u) &=& (u-q^{-3})(u-q^{-7}),\\
P_3(u) &=& u-q^{-10}, \\
P_4(u) &=& u - q^{-7},\\
P_k(u) &=& 1,\quad (5\le k\le N-1);
\end{eqnarray*} 
\begin{eqnarray*}
P_1(u) &=& (u-q^{-2})(u-q^{-4}),\\
P_2(u) &=& (u-q^{-7})(u-q^{-9}),\\
P_3(u) &=& P_4(u)\ \  =\ \ 1, \\
P_5(u) &=& u - q^{-8},\\
P_k(u) &=& 1,\quad (6\le k\le N-1);
\end{eqnarray*} 
\begin{eqnarray*}
P_1(u) &=& (u-q^{-2})(u-q^{-4}),\\
P_2(u) &=& u-q^{-7},\\
P_3(u) &=& u-q^{-10}, \\
P_4(u) &=& u - q^{-7},\\
P_k(u) &=& 1,\quad (5\le k\le N-1).
\end{eqnarray*} 
The modules $V(\nn_1)$ and $V(\nn_3)$ are non-zero 
only if $N\ge 6$. Hence $V(\mm_1)\otimes V(\mm_2)$
has only two composition factors $V(\nn_2)$ and
$V(\nn_4)$ for $N=5$, and four composition factors
$V(\nn_1), V(\nn_2)$, $V(\nn_3), V(\nn_4)$
for $N\ge 6$.
}
\end{example}

\subsection{}
We note that our result implies the following 
\begin{theorem}\label{th3}
All composition factors of the tensor product of two
evaluation modules of $U_q(\slchap_N)$ occur with multiplicity
one.
\end{theorem}


\bigskip
\small

\noindent
\begin{tabular}{ll}
{\sc B. Leclerc} : &
Laboratoire de Math\'ematiques Nicolas Oresme,\\
& Universit\'e de Caen, Campus II,\\
& Bld Mar\'echal Juin,
BP 5186, 14032 Caen cedex, France\\
&email : {\tt leclerc@math.unicaen.fr}
\end{tabular}

\end{document}